\documentclass[11pt]{article}
\usepackage{amsmath,amssymb}

\newtheorem{propo}{{\bf Proposition}}[section]
\newtheorem{coro}[propo]{{\bf Corollary}}
\newtheorem{lemma}[propo]{{\bf Lemma}} \newtheorem{theor}[propo]{{\bf
Theorem}} 

\def\Z{{\mathbb Z}}

\begin{document}

\vspace*{1.0in}

\begin{center} C-IDEALS OF LIE ALGEBRAS 
\end{center}
\bigskip

\begin{center} DAVID A. TOWERS 
\end{center}
\bigskip
\centerline {Department of
Mathematics, Lancaster University} \centerline {Lancaster LA1 4YF,
England} \centerline {Email: d.towers@lancaster.ac.uk}
\bigskip

\begin{abstract}
A subalgebra $B$ of a Lie algebra $L$ is called a {\em c-ideal} of $L$ if there is an ideal $C$ of $L$ such that $L = B + C$ and $B \cap C \leq B_L$, where $B_L$ is the largest ideal of $L$ contained in $B$. This is analogous to the concept of c-normal subgroup, which has been studied by a number of authors. We obtain some properties of c-ideals and use them to give some characterisations of solvable and supersolvable Lie algebras. We also classify those Lie algebras in which every one-dimensional subalgebra is a c-ideal.
\end{abstract}

\section{Introduction}
\medskip
Throughout $L$ will denote a finite-dimensional Lie algebra over a field $F$. 
If $B$ is a subalgebra of $L$ we define $B_L$, the {\em core} (with respect to $L$) of $B$ to be the largest ideal of $L$ contained in $B$. We say that a subalgebra $B$ of $L$ is a {\em c-ideal} of $L$ if there is an ideal $C$ of $L$ such that $L = B + C$ and $B \cap C \leq B_L$. This is analogous to the concept of c-normal subgroup as introduced by Wang in \cite{wang}; this concept has since been further studied by a number of authors, including Li and Guo (\cite{lg} and \cite{lg2}), Jehad (\cite{jehad}), Wang (\cite{wang2}), Wei (\cite{wei}) and Skiba (\cite{skiba}). 
\par
The maximal subalgebras of a Lie algebra $L$ and their relationship to the structure of $L$ have been studied extensively. It is known that $L$ is nilpotent if and only if every maximal subalgebra of $L$ is an ideal of $L$. A further result is that every maximal subalgebra of $L$ has codimension one in $L$ if and only if $L$ is supersolvable. In this paper we obtain some similar characterisations of solvable and supersolvable Lie algebras in terms of c-ideals.
\par

A subalgebra $B$ of $L$ is a {\em retract} of $L$ if there is an endomorphism $\theta : L \rightarrow L$ such that $\theta (b) = b$ for all $b \in B$ and $\theta (x) \in B$ for all $x \in L$. Such a map $\theta$ is called a {\em retraction}. Then it is easy to see that ideals of $L$ and retracts of $L$ are c-ideals of $L$; in the case of retracts the kernel of the retraction is an ideal that complements $B$. If $F$ has characteristic zero then every Levi factor of $L$ is a c-ideal of $L$.
\par
In section one we give some basic properties of c-ideals; in particular, it is shown that c-ideals inside the Frattini subalgebra of a Lie algebra $L$ are necessarily ideals of $L$. In section two we first show that all maximal subalgebras of $L$ are c-ideals of $L$ if and only if $L$ is solvable. It is further shown that, over a field of characteristic zero or over an algebraically closed field of characteristic $p >5$, $L$ has a solvable maximal subalgebra that is a c-ideal if and only if $L$ is solvable. Finally we have that if all maximal nilpotent subalgebras of $L$ are c-ideals, or if all Cartan subalgebras of $L$ are c-ideals and $F$ has characteristic zero, then $L$ is solvable.
\par
In section three we show that if every maximal subalgebra of each maximal nilpotent subalgebra of $L$ is a c-ideal of $L$ then $L$ is supersolvable. If each of the maximal nilpotent subalgebras of $L$ has dimension at least two then the assumption of solvability can be removed. Similarly if the field has characteristic zero and $L$ is not three-dimensional simple then this restriction can be removed. In the final section we classify those Lie algebras in which every one-dimensional subalgebra is a c-ideal. 
\par

If $A$ and $B$ are subalgebras of $L$ for which $L = A + B$ and $A \cap B = 0$ we will write $L = A \oplus B$. The ideals $L^{(k)}$ and $L^k$ are defined inductively by $L^{(1)} = L^1 = L$, $L^{(k+1)} = [L^{(k)},L^{(k)}]$, $L^{k+1} = [L,L^k]$ for $k \geq 1$. If $A$ is a subalgebra of $L$, the {\em centralizer} of $A$ in $L$ is $C_{L}(A) = \{ x \in L : [x, A] = 0\}$.  
\bigskip

\section{Preliminary results}
\medskip
First we give some basic properties of c-ideals.
\bigskip

\begin{lemma}\label{l:ci}
\begin{itemize}
\item[(i)] If $B$ is a c-ideal of $L$ and $B \leq K \leq L$ then $B$ is a c-ideal of $K$.
\item[(ii)] If $I$ is an ideal of $L$ and $I \leq B$ then $B$ is a c-ideal of $L$ if and only if $B/I$ is a c-ideal of $L/I$.
\end{itemize}
\end{lemma}
\medskip
{\it Proof.}
\begin{itemize}
\item[(i)] Suppose that $B$ is a c-ideal of $L$ and $B \leq K \leq L$. Then there is an ideal $C$ of $L$ with $L = B + C$ and $B \cap C \leq B_L$. It follows that $K = (B + C) \cap K = B + C \cap K$, where $C \cap K$ is an ideal of $K$ and $B \cap C \cap K \leq B_L \cap K \leq B_K$, and so $B$ is a c-ideal $K$.
\item[(ii)] Suppose first that $B/I$ is a c-ideal of $L/I$. Then there is an ideal $C/I$ of $L/I$ such that $L/I = B/I + C/I$ and $(B/I) \cap (C/I) \leq (B/I)_{L/I} = B_L/I$. It follows that $L = B + C$, where $C$ is an ideal of $L$ and $B \cap C \leq B_L$, whence $B$ is a c-ideal of $L$.
\par
Suppose conversely that $I$ is an ideal of $L$ with $I \leq B$ such that $B$ is a c-ideal of $L$. Then there is an ideal $C$ of $L$ such that $L = B + C$ and $B \cap C \leq B_L$. Now $L/I = B/I + (C + I)/I$, where $(C + I)/I$ is an ideal of $L/I$ and
$(B/I) \cap (C + I)/I = (B \cap (C + I))/I = (I + B \cap C)/I \leq B_L/I = (B/I)_{L/I}$, so $B/I$ is a c-ideal of $L/I$.
\end{itemize} 
\bigskip

The {\em Frattini subalgebra} of $L$, $F(L)$, is the intersection of all of the maximal subalgebras of $L$. The {\em Frattini ideal}, $\phi(L)$, of $L$ is $F(L)_L$.  The next result shows, in particular, that c-ideals inside the Frattini subalgebra of a Lie algebra $L$ are necessarily ideals of $L$.
\bigskip

\begin{propo}\label{p:frat}
Let $B, C$ be subalgebras of $L$ with $B \leq F(C)$. If $B$ is a c-ideal of $L$ then $B$ is an ideal of $L$ and $B \leq \phi(L)$.
\end{propo}
\medskip
{\it Proof.} Suppose that $L = B + K$ and $B \cap K \leq B_L$. Then $C = C \cap L = C \cap (B + K) = B + C \cap K = C \cap K$ since $B \leq F(C)$. Hence $B \leq C \leq K$, giving $B = B \cap K \leq B_L$ and $B$ is an ideal of $L$. It then follows from \cite[Lemma 4.1]{frat} that $B \leq \phi(L)$.
\bigskip

\section{Some characterisations of solvable algebras}
\medskip

\begin{theor}\label{t:max}
Let $L$ be a Lie algebra over any field $F$. Then all maximal subalgebras of $L$ are c-ideals of $L$ if and only if $L$ is solvable.
\end{theor}
\medskip
{\it Proof.} Let $L$ be a non-solvable Lie algebra of smallest dimension in which maximal subalgebras are c-ideals of $L$. Then all proper factor algebras of $L$ are solvable, by Lemma \ref{l:ci} (ii). Suppose first that $L$ is simple. Let $M$ be a maximal subalgebra of $L$. Then $M$ is a c-ideal so there is an ideal $C$ of $L$ such that $L = M + C$ and $M \cap C \leq M_L = 0$, as $L$ is simple. This yields that $C$ is a non-trivial proper ideal of $L$, a contradiction. If $L$ has two minimal ideals $B_1$ and $B_2$, then $L/B_1$ and $L/B_2$ are solvable and $B_1 \cap B_2 = 0$, so $L$ is solvable. Hence $L$ has a unique minimal ideal $B$ and $L/B$ is solvable.
\par
Suppose there is an element $b \in B$ such that ${\rm ad}_L b$ is not nilpotent. Let $L = L_0 \oplus L_1$ be the Fitting decomposition of $L$ relative to ${\rm ad}_L b$. Then $L \neq L_0$ so let $M$ be a maximal subalgebra of $L$ containing $L_0$. As $M$ is a c-ideal there is an ideal $C$ of $L$ such that $L = M + C$ and $M \cap C \leq M_L$. Now $L_1 \leq B$ so $B \not \leq M_L$. It follows that $M_L = 0$ whence $M = L_0$ and $B = C = L_1$. But $b \in M \cap B = 0$. Hence every element of $B$ is ad-nilpotent, yielding that $B$ is nilpotent and so $L$ is solvable, a contradiction. 
\par
Now suppose that $L$ is solvable and let $M$ be a maximal subalgebra of $L$. Then there is a $k \geq 2$ such that $L^{(k)} \leq M$, but $L^{(k-1)} \not \leq M$. We have that $L^{(k-1)}$ is an ideal of $L$, $L = M + L^{(k-1)}$ and $M \cap L^{(k-1)} \leq M_L$, so $M$ is a c-ideal of $L$.
\bigskip

\begin{theor}\label{t:solv}
Let $L$ be a Lie algebra over a field $F$ of characteristic zero. Then $L$ has a solvable maximal subalgebra that is a c-ideal of $L$ if and only if $L$ is solvable.
\end{theor}
\medskip
{\it Proof.} Suppose first that $L$ has a solvable maximal subalgebra $M$ that is a c-ideal of $L$. We show that $L$ is solvable. Let $L$ be a minimal counter-example. Then there is an ideal $K$ of $L$ such that $L = M + K$ and $M \cap K \leq M_L$. Now $M_L = 0$, since otherwise, $L/M_L$ is solvable and $M_L$ is solvable, whence $L$ is solvable, a contradiction. It follows that $L = M \oplus K$. If $R$ is the solvable radical of $L$ then $R \leq M_L = 0$, so $L$ is semisimple and $L^2 = L$. But $L^2 \leq M^2 + K \neq L$, a contradiction. The result follows.
\par
The converse follows from Theorem \ref{t:max}.
\bigskip

For fields of characteristic $p > 0$ we have the following result.
\bigskip

\begin{theor}\label{t:solvp}
Let $L$ be a Lie algebra over an algebraically closed field $F$ of characteristic greater than 5. Then $L$ has a solvable maximal subalgebra that is a c-ideal of $L$ if and only if $L$ is solvable.
\end{theor}
\medskip
{\it Proof.} Suppose first that $L$ has a solvable maximal subalgebra $M$ that is a c-ideal of $L$. We show that $L$ is solvable. Let $L$ be a minimal counter-example. Then, as above, $L = M \oplus K$ and $K$ is a minimal ideal of $L$. We follow the contents of \cite{weis}: $M$ defines a filtration in which $L_0 = M$, $L_{i+1} = \{x \in L_i : [x, L] = L_i \}$. When $L_1 = 0$ this filtration is called {\em short}; otherwise it is {\em long}. Suppose first that it is short. Then, as in the first two paragraphs of the proof of \cite[Theorem 2.2]{weis}, $L = \bigoplus_{i \in \Z_p} L_i$, $M = L_0$ and $L_i$ is an irreducible $M$-submodule of $L$ for each $i \neq 0$. Moreover, since $K$ is an ideal of $L$, $K = \bigoplus_{0 \neq i \in \Z_p} L_i$. Let $S$ be the subalgebra generated by $L_1$. Then $S$ is spanned by commutators $c(x_1, \ldots, x_n) = [x_1,[x_2, \ldots [x_{n-1}, x_n]]$ with $x_i \in L_1$ and $n \geq 1$. Now $c(x_1, \ldots, x_p) \in M \cap K = 0$ for all $x_1, \ldots, x_p \in L_1$, so $S$ is nilpotent. Also $M$ idealizes $S$, so $M + S$ is a subalgebra of $L$, whence $L = M + S$ and $S$ is an ideal of $L$. It follows that $K = S$ is nilpotent and $L$ is solvable, a contradiction.
\par
Now suppose that the filtration is long. Then the nilradical, $N$, of $M$ acts nilpotently on $K$, by \cite[Proposition 2.5]{weis}. Let $C = C_{K}(N)$. Then $C \neq 0$ and $M + C$ is a subalgebra of $L$. It follows that $L = M + C$, whence $K = C$..  But this means that $N$ is an ideal of $L$, so that $N \subseteq M_L = 0$. We conclude that $M = 0$, a contradiction.
\par
The converse follows from Theorem \ref{t:max} as before.
\bigskip

\begin{theor}\label{t:maxnilp}
Let $L$ be a Lie algebra over any field $F$, such that all maximal nilpotent subalgebras of $L$ are c-ideals of $L$. Then $L$ is solvable.
\end{theor}
\medskip
{\it Proof.} Let $N$ be the nilradical of $L$ and let $x \notin N$. Then $x \in B$ for some maximal nilpotent subalgebra $B$ of $L$, and there is an ideal $C$ of $L$ such that $L = B + C$ and $B \cap C \leq B_L$. Clearly $x \notin B_L \leq N$, so $x \notin C$. Moreover, $L/C \cong B/(B \cap C)$ is nilpotent. So if $x \notin N$, there is an ideal $C$ of $L$ such that $x \notin C$ and $L/C$ is nilpotent.
\par
So let $x_1 \notin N$ and let $C_1$ be such an ideal with $x_1 \notin C_1$ and $L/C_1$ nilpotent. If $C_1 \leq N$ we have finished. If not, then choose $x_2 \in C_1 \setminus N$ and let $C_2$ be such an ideal with $x_2 \notin C_2$ and $L/C_2$ nilpotent. Clearly dim $(C_1 \cap C_2) <$ dim $C_1$. If $C_1 \cap C_2 \not \leq N$, choose $x_3 \in (C_1 \cap C_2) \setminus N$. Continuing in this way we find ideals $C_1, \hdots , C_n$ of $L$ such that $C_1 \cap \hdots \cap C_n \leq N$ and $L/C_i$ is nilpotent for each $1 \leq i \leq n$. Since $L/(C_1 \cap \hdots \cap C_n)$ is solvable, the result follows. 
\bigskip

\begin{theor}\label{t:cart}
Let $L$ be a Lie algebra, over a field $F$ of characteristic zero, in which every Cartan subalgebra of $L$ is a c-ideal of $L$. Then $L$ is solvable.
\end{theor}
\medskip
{\it Proof.} Suppose that every Cartan subalgebra of $L$ is a c-ideal of $L$, and that $L$ has a non-zero Levi factor $S$. Let $H$ be a Cartan subalgebra of $S$ and let $B$ be a Cartan subalgebra of its centralizer in the solvable radical of $L$. Then $C = H + B$ is a Cartan subalgebra of $L$ (see \cite{dix}) and  there is an ideal $K$ of $L$ such that $L = C + K$ and $C \cap K \leq C_L$. Now there is an $r \geq 2$ such that $L^{(r)} \leq K$. But $S \leq L^{(r)} \leq K$, so $C \cap S \leq C \cap K \leq C_L$ giving $C \cap S \leq C_L \cap S = 0$, a contradiction. It follows that $S = 0$ and hence that $L$ is solvable.   
\bigskip

\noindent {\bf Note}: If $L^{\infty} = \cap_{i=1}^{\infty} L^i$ is abelian then the converse to the above theorem holds, by Theorem 4.4.1.1 of \cite{wint}.
\bigskip

\section{Some characterisations of supersolvable algebras}
\medskip
First we need some preliminary results concerning maximal nilpotent subalgebras of Lie algebras.
\bigskip

\begin{lemma}\label{l:maxnilp}
Let $L$ be a Lie algebra over any field $F$, let $A$ be an ideal of $L$ and let $U/A$ be a maximal nilpotent subalgebra of $L/A$. Then $U = C+A$, where $C$ is a maximal nilpotent subalgebra of $L$.
\end{lemma}
\medskip
{\it Proof.} If $A \leq \phi(U)$ then $U/ \phi(U)$ is nilpotent, whence $U$ is nilpotent, by Theorem 6.1 of \cite{frat} and the result is clear. So suppose that $A \not \leq \phi(U)$. Then $U = A + M$ for some maximal subalgebra $M$ of $U$. If we choose $B$ to be minimal with respect to $U = A + B$, then $A \cap B \leq \phi(B)$, by Lemma 7.1 of \cite{frat}. Also $U/A \cong B/(A \cap B)$ is nilpotent, which yields that $B$ is nilpotent. If we now choose $C$ to be the biggest nilpotent subalgebra of $U$ such that $U = A + C$, it is easy to see that $C$ is a maximal nilpotent subalgebra of $L$.   
\bigskip

\begin{lemma}\label{l:reduce}
Let $L$ be a Lie algebra, over any field $F$, in which every maximal subalgebra of each maximal nilpotent subalgebra of $L$ is a c-ideal of $L$, and let $A$ be a minimal abelian ideal of $L$. Then every maximal subalgebra of each maximal nilpotent subalgebra of $L/A$ is a c-ideal of $L/A$.
\end{lemma}
\medskip
{\it Proof.} Suppose that $U/A$ is a maximal nilpotent subalgebra of $L/A$. Then $U = C + A$ where $C$ is a maximal nilpotent subalgebra of $L$ by Lemma \ref{l:maxnilp}. Let $B/A$ be a maximal subalgebra of $U/A$. Then $B = B \cap (C+A) = B \cap C + A = D + A$ where $D$ is a maximal subalgebra of $C$ with $B \cap C \leq D$. Now $D$ is a c-ideal of $L$ so there is an ideal $K$ of $L$ with $L = D+K$ and $D \cap K \leq D_L$. 
\par

If $A \leq K$ we have $L/A = (D+K)/A = ((D+A)/A) + (K/A) = (B/A) + (K/A)$, and $(B/A) \cap (K/A) = (B \cap K)/A = ((D+A) \cap K)/A = (D \cap K + A)/A \leq (D_L + A)/A \leq (B/A)_{L/A}$. 
\par

If $A \not \leq K$, we have $A \cap K = 0$. Then $(A+K)/K$ is a minimal ideal of $L/K$, which is nilpotent, so dim$A = 1$ and $LA \leq A \cap K = 0$.  It follows that $A \leq C$ and $B = D$. We have $L = B+K$ and $B \cap K \leq B_L$, so $L/A = (B/A) + ((K+A)/A)$ and $(B/A) \cap ((K+A)/A) = (B \cap (K+A))/A = (B \cap K + A)/A \leq (B_L + A)/A \leq (B/A)_{L/A}$.
\bigskip

\begin{lemma}\label{l:onedim}
Let $L$ be a Lie algebra, over any field $F$, in which every maximal subalgebra of each maximal nilpotent subalgebra of $L$ is a c-ideal of $L$, and suppose that $A$ is a minimal abelian ideal of $L$ and $M$ is a core-free maximal subalgebra of $L$. Then $A$ is one dimensional.
\end{lemma}
\medskip
{\it Proof.} We have that $L = A \oplus M$ and $A$ is the unique minimal ideal of $L$, by Lemma 1.4 of \cite{tow}. Let $C$ be a maximal nilpotent subalgebra of $L$ with $A \leq C$. If $C = A$, choose $B$ to be a maximal subalgebra of $A$, so that $A = B + Fa$ and $B_L = 0$. Then $B$ is a c-ideal of $L$ so there is an ideal $K$ of $L$ with $L = B + K$ and $B \cap K \leq B_L = 0$. But now $L = A + K = K$, giving $B = 0$ and dim$A = 1$.
\par

So suppose that $C \neq A$. Then $C = A + M \cap C$. Let $B$ be a maximal subalgebra of $C$ containing $M \cap C$. Then $B$ is a c-ideal of $L$, so there is an ideal $K$ of $L$ with $L = B + K$ and $B \cap K \leq B_L$. If $A \leq B_L \leq B$ we have $C = A + M \cap C \leq B$, a contradiction. Hence $B_L = 0$ and $L = B \oplus K$. Now $C = B + C \cap K$ and $B \cap C \cap K = B \cap K = 0$. As $C$ is nilpotent this means that dim$(C \cap K) = 1$. But $A \leq C \cap K$, so dim$A = 1$, as required.
\bigskip

We can now prove our main result.
\bigskip

\begin{theor}\label{t:sup}
Let $L$ be a solvable Lie algebra, over any field $F$, in which every maximal subalgebra of each maximal nilpotent subalgebra of $L$ is a c-ideal of $L$. Then $L$ is supersolvable.
\end{theor}
\medskip
{\it Proof.}  Let $L$ be a minimal counter-example and let $A$ be a minimal abelian ideal of $L$. Then $L/A$ satisfies the same hypotheses, by Lemma \ref{l:reduce}. We thus have that $L/A$ is supersolvable, and it remains to show that dim$A = 1$.
\par
If there is another minimal ideal $I$ of $L$, then $A \cong (A+I)/I \leq L/I$, which is supersolvable and so dim$A = 1$. So we can assume that $A$ is the unique minimal ideal of $L$. Also, if $A \leq \phi(L)$ we have that $L/ \phi(L)$ is supersolvable, whence $L$ is supersolvable, by Theorem 7 of \cite{barnes2}. We therefore further assume that $A \not \leq \phi(L)$. It follows that $L = A \oplus M$, where $M$ is a core-free maximal subalgebra of $L$.The result now follows from Lemma \ref{l:onedim}.
\bigskip

If $L$ has no one-dimensional maximal nilpotent subalgebras, we can remove the solvability assumption from the above result.
\bigskip

\begin{coro}\label{c:sup}
Let $L$ be a Lie algebra, over any field $F$, in which every maximal nilpotent subalgebra has dimension at least two. If every maximal subalgebra of each maximal nilpotent subalgebra of $L$ is a c-ideal of $L$, then $L$ is supersolvable.
\end{coro}
\medskip
{\it Proof.} Let $N$ be the nilradical of $L$ and let $x \notin N$. Then $x \in C$ for some maximal nilpotent subalgebra $C$ of $L$. Since dim$ C > 1$, there is a maximal subalgebra $B$ of $C$ with $x \in B$. Now there is an ideal $K$ of $L$ with $L = B + K$ and $B \cap K \leq B_L \leq C_L \leq N$. Clearly $x \notin K$, since otherwise $x \in B \cap K \leq N$. Moreover, $L/K$ is nilpotent. We have shown that if $x \notin N$ there is an ideal $K$ of $L$ with $x \notin K$ and $L/K$ nilpotent. Proceeding as in Theorem \ref{t:maxnilp} we see that $L$ is solvable. The result then follows from Theorem \ref{t:sup}.   
\bigskip

If $L$ has a one-dimensional maximal nilpotent subalgebra then we can also remove the solvability assumption from Theorem \ref{t:sup} provided that the underlying field $F$ has characteristic zero and $L$ is not three-dimensional simple.
\bigskip

\begin{coro}\label{c:zero}
Let $L$ be a Lie algebra over a field $F$ of characteristic zero. If every maximal subalgebra of each maximal nilpotent subalgebra of $L$ is a c-ideal of $L$, then $L$ is supersolvable or three-dimensional simple.
\end{coro}
\medskip
{\it Proof.} If every maximal nilpotent subalgebra of $L$ has dimension at least two then $L$ is supersolvable, by Corollary \ref{c:sup}. So we need only consider the case where $L$ has a one-dimensional maximal nilpotent subalgebra, $Fx$ say. 
\par
Suppose first that $L$ is semisimple, so $L = S_1 \oplus \ldots \oplus S_n$, where $S_i$ is a simple ideal of $L$ for $1 \leq i \leq n$. Let $n > 1$. If $x \in S_i$ then choosing $s \in S_j$ with $j \neq i$ we have that $Fx + Fs$ is a two-dimensional abelian subalgebra, which contradicts the maximality of $Fx$. If $x \notin S_i$ for every $1 \leq i \leq n$, then $x$ has non-zero projections in at least two of the $S_k$'s, say $s_i \in S_i$ and $s_j \in S_j$. But then $Fx + Fs_i$ is a two-dimensional abelian subalgebra, a contradiction again. It follows that $L$ is simple. But then $Fx$ is a Cartan subalgebra of $L$, which yields that $L$ has rank one and thus is three dimensional.
\par
So now let $L$ be a minimal counter-example. We have seen that $L$ is not semisimple, so it has a minimal abelian ideal $A$. By Lemma \ref{l:reduce}, $L/A$ is supersolvable or three-dimensional simple. In the former case, $L$ is solvable and so supersolvable, by Theorem \ref{t:sup}. In the latter case, $L = A \oplus S$ where $S$ is three-dimensional simple, and so a core-free maximal subalgebra of $L$. It follows from Lemma \ref{l:onedim} that dim$A = 1$. But now $C_L(A) = A$ or $L$. In the former case $S \cong L/A = L/C_L(A) \cong Inn(A)$, a subalgebra of $Der(A)$, which is impossible. Hence $L = A \oplus S$, where $A$ and $S$ are both ideals of $L$, and again $L$ has no one-dimensional maximal nilpotent subalgebras.   
\bigskip

\section{One-dimensional c-ideals}
\medskip
First we note that one-dimensional c-ideals are easy to classify.
\bigskip

\begin{propo}\label{p:min}
Let $L$ be a Lie algebra over any field $F$. Then the one-dimensional subalgebra $Fx$ of $L$ is a c-ideal of $L$ if and only if 
\begin{itemize}
\item[(i)] $Fx$ is an ideal of $L$; or
\item[(ii)] $x \notin L^2$.
\end{itemize}
\end{propo}
\medskip
{\it Proof.} Let $Fx$ be a c-ideal of $L$. Then there is an ideal $K$ of $L$ such that $L = Fx + K$ and $Fx \cap K \leq (Fx)_L$. But $Fx \cap K = Fx$ or $0$. The former implies that $Fx$ is an ideal of $L$, and the latter implies that $x \notin L^2 \leq K$.
\par

Conversely, suppose that $x \notin L^2$. Then there is a subspace $K$ of $L$ of codimension one in $L$ such that $L^2 \leq K$ and $x \notin K$. Clearly $L = Fx \oplus K$ and $K$ is an ideal of $L$, whence $Fx$ is a c-ideal of $L$.
\bigskip 

We shall denote by $Z(L)$ the {\em centre} of $L$; that is $Z(L) = \{x \in L : [x,y] = 0 \hbox{ for all } y \in L \}$. The {\em abelian socle} of $L$, Asoc$L$, is the sum of the minimal abelian ideals of $L$. We say that $L$ is {\em almost abelian} if $L = L^2 \oplus Fx$, where $L^2$ is abelian and $[x,y] = y$ for all $y \in L^2$.
\bigskip

\begin{theor}\label{t:min}
Let $L$ be a Lie algebra over any field $F$. Then all one-dimensional subalgebras of $L$ are c-ideals of $L$ if and only if 
\begin{itemize}
\item[(i)] $L^3 = 0$; or
\item[(ii)] $L = A \oplus B$, where $A$ is an abelian ideal of $L$ and $B$ is an almost abelian ideal of $L$.
\end{itemize}
\end{theor}
\medskip
{\it Proof.} Suppose that all one-dimensional subalgebras of $L$ are c-ideals of $L$. First note that the one-dimensional ideals are inside Asoc$L$. If $Fx$ is not an ideal of $L$ then there is an ideal $M$ of $L$ such that $L = Fx + M$ and $Fx \cap M \leq (Fx)_L = 0$, so $Fx$ is complemented by an ideal of codimension one in $L$. 
\par
Now let $A$ be a minimal ideal of $L$ and let $a \in A$. If $A \neq Fa$ then there is an ideal $M$ of codimension one in $L$ which complements $Fa$. But this implies that $M \cap A = 0$, whence $A = Fa$, a contradiction. It follows that all minimal ideals are one dimensional. Put Asoc$L = Fa_1 \oplus \ldots \oplus Fa_r$. Suppose that $[x,a_i] = \lambda a_i$, $[x,a_j] = \mu a_j$ and $\lambda \neq \mu$. Then $F(a_i + a_j)$ is not an ideal of $L$, and so there is an ideal $M$ of $L$ with $L = F(a_i + a_j) \oplus M$. Clearly one of $a_i, a_j$ does not belong to $M$: suppose $a_i \notin M$. Then $L = Fa_i \oplus M$ and $a_i \in Z(L)$. Hence Asoc$L = Z \oplus D$, where $Z = Z(L)$ and $[x,a] = \lambda_x a$ for all $a \in D$ and $\lambda_x \neq 0$.
\par
Let $\Lambda : L \rightarrow F$ be given by $\Lambda (x) = \lambda_x$. This is a one-dimensional representation of $L$. Hence, either Im $\Lambda = 0$, in which case $D = 0$, or else $L =$ Ker $\Lambda \oplus Fx$ and $\lambda_x = 1$. Put $L = Z \oplus D \oplus C \oplus Fx$, where $C \subseteq$ Ker $\Lambda$.
\par
If $y \notin \hbox{Asoc} L$ then $Fy$ is complemented by an ideal $M$ and $L^2 \leq M$, so $y \notin L^2$. This yields that $L^2 \leq \hbox{Asoc} L$. Clearly $D \leq L^2$. If $L^2 \leq Z$ then $L^3 = 0$ and we have case (i). So suppose that $D \neq 0$ and let $a \in D$. If there is an element $z \in Z \cap L^2$, then $F(z+a)$ is not an ideal of $L$ and so $z+a \notin L^2$, a contradiction. Hence $L^2 = D$.
\par
Let $c \in C$. If $[x,c] = 0$ then $Fc$ is an ideal of $L$ and $c \in C \cap \hbox{Asoc} L = 0$. So suppose $[x,c] \neq 0$. Then $[x,c] \in D$ so $[x,c-[x,c]] = 0$. This implies that $F(c - [x,c])$ is an ideal of $L$, whence $c - [x,c] \in  \hbox{Asoc} L$. But now $c \in C \cap \hbox{Asoc} L = 0$. Hence $C = 0$ and $L$ is as described in (ii). 
\par 
Now suppose that $L^3 = 0$. If $x \in L^2$ then $Fx$ is an ideal of $L$. If $x \notin L^2$ there is a subspace $M$ of codimension one in $L$ containing $L^2$ such that $x \notin M$. This implies that $Fx$ is a c-ideal of $L$.
\par
Finally suppose that $L$ is as in (ii): say $L = Z \oplus A \oplus Fx$ where $Z = Z(L)$, $A$ is abelian and $[x,a] = a$ for all $a \in A$. Let $z+a+ \alpha x \in L$. If $z \neq 0$ then choosing $M = Z_1 \oplus A \oplus Fx$ where $Z = Z_1 \oplus Fz$ shows that $F(z+a+ \alpha x)$ is a c-ideal of $L$. So suppose $z = 0$. If $\alpha = 0$ then $Fa$ is an ideal of $L$. If $\alpha \neq 0$ then choosing $M = Z \oplus A$ shows that $F(a+ \alpha x)$ is a c-ideal of $L$.

\end{document}